\documentclass[a4paper]{amsart}
 
\usepackage{amscd,amssymb}

\newcommand{\bb}{\mathbb}
\newcommand{\cc}{\bb C}
\newcommand{\rr}{\bb R}
\newcommand{\z}{\bb Z}
\newcommand{\zz}{\z/2}
\newcommand{\pp}{\bb P}
\newcommand{\kk}{\bb K}

\renewcommand{\emptyset}{\varnothing}
 
\newcommand{\labelto}[2]{\xrightarrow{\makebox[#2 em]{\scriptsize
${#1}$}}}
 
\newcommand{\xsig}{X}
 
\newcommand{\xr}{X(\rr)}
 
\newcommand{\xc}{X(\cc)}

\newcommand{\xcomp}{X_\cc}
 
\newcommand{\yc}{Y(\cc)}

\newcommand{\yr}{Y(\rr)}

\newcommand{\cL}{\mathcal{L}}
\newcommand{\cE}{\mathcal{E}}

\newcommand{\pic}{\operatorname{Pic}}

\newcommand{\ns}{\operatorname{NS}}

\newcommand{\im}{\operatorname{Im}}
\newcommand{\kod}{\operatorname{kod}}
\newcommand{\elm}{\operatorname{elm}}
 
\newcommand{\ox}{\mathcal O_X}
\newcommand{\kx}{\mathcal K_X}

\newcommand{\halg}{H_{\mathrm{alg}}}
\newcommand{\hcalg}{H_{\cc-\mathrm{alg}}}

\newcommand{\refbooktitle}[1]{\textit{#1}}
\newcommand{\refjourtitle}[1]{\textit{#1}}
\newcommand{\refpapertitle}{}
 
\numberwithin{equation}{section}
\newtheorem{theo}[equation]{Theorem}%
\newtheorem{prop}[equation]{Proposition}
\newtheorem{lem}[equation]{Lemma}
\newtheorem{cor}[equation]{Corollary}

\newtheorem{defin}[equation]{Definition}
\theoremstyle{remark}
\newtheorem{rem}[equation]{Remark}

\begin{document}
 
\title[Real conic bundles]{Real algebraic morphisms on 2-dimensional conic bundles}
 
\author{Fr\'ed\'eric Mangolte}
\thanks{I would like to thank J.~Bochnak for very helpful discussions. 
Many thanks also to C.~Raffalli and P.~Verovic for their careful reading.}


\address{Fr\'ed\'eric Mangolte, Laboratoire de Math\'ematiques,
Universit\'e de Savoie, 73376 Le Bourget-du-Lac Cedex, France,
T\'el~: (33) 4 79 75 86 60, Fax~: (33) 4 79 75 81 42}
\email{mangolte@univ-savoie.fr}

\begin{abstract}
Given two nonsingular real algebraic
varieties $V$ and $W$, we consider the problem of deciding whether a smooth map $f:V\rightarrow W$ can be
approximated by regular maps in the space of
${\mathcal{C}}^\infty$ mappings from $V$ to $W$ in the
${\mathcal{C}}^\infty$ topology.
 
Our main result is a complete solution to this problem in case $W$ is the usual 2-dimensional sphere and $V$ is a real algebraic surface of negative Kodaira dimension. 
\end{abstract}
 
\maketitle
 
\section{Introduction}\label{sec:intro}
 
We deal with regular maps between real algebraic varieties. Here the term {\em real algebraic variety} stands for a locally ringed space isomorphic to a locally closed subset of $\pp^n(\rr)$, for some $n$, endowed with
the Zariski topology and the sheaf of $\rr$-valued regular
functions. Morphisms between real algebraic varieties are called
{\em regular maps}. 

An equivalent description of real algebraic varieties can be obtained using quasi-projective varieties defined over $\rr$. Given such a variety $X$, the Galois
group $G = \{1, \sigma\}$ of $\cc \vert \rr$ acts on $\xc$, the set of complex
points of $X$, via an antiholomorphic involution.
The real part $\xr$ is  then
precisely the set of fixed points under this action. If $\xr$ is Zariski dense in $X$, then we consider it as a real algebraic variety whose structure sheaf is the restriction of the structure sheaf of $X$. Therefore, a regular map $V\to W$ in the above sense is the restriction of a rational map $\xc\dasharrow\yc$ with no poles on $V=\xr$.

Every real algebraic variety is isomorphic to
a Zariski closed subvariety of $\rr ^n$. Thus, all topological notions
about real algebraic varieties will refer to the Euclidean
topology of $\rr ^n$.

Given two nonsingular real algebraic varieties $V$ and $W$, with
$V$ compact, we consider the set $\mathcal {R}(V,W)$ of all
regular maps from $V$ into $W$ as a subset of the space $\mathcal
{C}^\infty (V,W)$ of all $\mathcal {C}^\infty$ maps from $V$ into
$W$ equipped with the $\mathcal {C}^\infty$ topology. We want to
study which $\mathcal {C}^\infty$ maps from $V$ into $W$ can be
approximated by regular maps. The classical Stone-Weierstrass
approximation theorem implies that all $\mathcal {C}^\infty$ maps
from $V$ into $W$ can be approximated by regular maps when $W=\rr
^m $ for some $m$. In this paper, we will mainly consider the case when
$W=S^2 $ is the usual two dimensional Euclidean sphere.

\bigskip
The main achievement of this paper is a complete answer to the approximation problem of smooth maps from a real algebraic surface of negative Kodaira dimension  into the 2-sphere.

In dimension less than three, an algebraic variety of negative Kodaira dimension is $\cc$-uniruled. In dimension two, such a variety is $\cc$-ruled. By definition, a $\cc$-ruled surface $X$ is $\cc$-birationally equivalent to a product $\pp^1\times B$, where $B$ is a complex  algebraic curve. When the genus $g(B)$ of $B$ is non vanishing, the surface $X$ admits a relatively minimal model $Y$ over $\cc$ endowed with a $\pp^1$-bundle structure $Y\longrightarrow B$.

When $X$ is defined over $\rr$, $X$ may not be $\rr$-birational to a product though it is $\cc$-birational to $\pp^1\times B$; consider for example a maximal real Del Pezzo surface of degree 2 or 1. 

Even when $g(B)\ne 0$, it may occur that no relatively $\rr$-minimal model $Y$ of $X$ is a $\pp^1$-bundle.  In this case, the surface $X$ admits only relatively $\rr$-minimal models which are not $\cc$-minimal: these are the real  conic bundles, see Section~\ref{sec:conic}.

\bigskip
Among all the surfaces of negative Kodaira dimension are the rational surfaces. There are several ways to define real algebraic varieties, hence several
ways to define real rational surfaces. A real algebraic surface $V$ is {\em $\cc$-rational} (or geometrically rational) if its complexification $X$ is $\cc$-birationally equivalent to $\pp^2_\cc$. Similarly, the surface $V$ is
{\em $\rr$-rational} if $V$ is $\rr$-birationally equivalent to $\pp^2_\rr$.

The approximation problem of smooth maps from $\rr$-rational surfaces into the 2-sphere was solved by J.~Bochnak and
W.~Kucharz \cite{BK87a,BK87b}.

In collaboration with N.~Joglar, we generalized this to the case when the source space is a $\cc$-rational real algebraic surface \cite{JM03}.

In this paper, we solve the problem for all the surfaces of negative Kodaira dimension.

Let $X$ be a real algebraic surface of negative Kodaira dimension which is not $\cc$-rational. Hence $X$ admits a real ruling $\rho \colon X \longrightarrow B$. Recall that a connected component of $\xr$ may be diffeomorphic to a torus, a sphere or any nonorientable surface. We denote by $K'$ the (possibly empty) set of connected components of $\xr$ which are diffeomorphic to the Klein bottle and whose image by $\rho$ is a connected component of $B(\rr)$.

\begin{theo}
\label{theo:mainf}
Let $X$ be a $\cc$-ruled non $\cc$-rational real algebraic surface. Given a smooth map $f\colon \xr \longrightarrow S^2$, the following conditions are equivalent:

\begin{enumerate}
\item $f$ can be approximated by regular maps;
\item $f$ is homotopic to a regular map;
\item for each component $M$ of $\xr$ diffeomorphic to a torus, $\deg(f)_{\vert M}=0$  and for each pair of components belonging to $K'$, $\deg_{\zz}(f)_{\vert M}=\deg_{\zz}(f)_{\vert N}$ .
\end{enumerate}

\end{theo}

In \cite{Ku}, W.~Kucharz gave another kind of generalization of his result with J.~Bochnak about $\rr$-rational surfaces. Namely, he solved the approximation problem of smooth maps from $\rr$-rational surfaces into $\rr$-rational surfaces.

We extend this result as follows (recall that a $\cc$-rational real algebraic surface $X$  is  $\rr$-rational if and only if $\xr$ is connected):

\begin{theo}
\label{theo:ruledrat}
Let $V=\xr$ and $W=\yr$  be  connected real algebraic surfaces such that $X$ is $\cc$-ruled and $Y$ is $\cc$-rational. Then the 
space $\mathcal {R}(V,W)$ is dense in the space $\mathcal {C}^\infty(V,W)$, except when $V$ is diffeomorphic to a torus and $W$ is diffeomorphic to a sphere.

In the latter case, the closure of $\mathcal {R}(V,W)$ in $\mathcal {C}^\infty(V,W)$ consists precisely of the null homotopic maps.
\end{theo}

Furthermore, we have answered a question raised by J. Bochnak. The following result is largely independent of the previous ones.

\begin{theo}
\label{theo:klein}
Let $V$ be a $\rr$-rational real algebraic surface diffeomorphic to the Klein bottle. Then $V$ is biregularly isomorphic to the blow-up of the real projective plane over one point. In other words, a $\rr$-rational Klein surface always admits  a non minimal smooth complexification.
\end{theo}

This result fits with the more general setting: a rational model of  a connected compact variety $M$ is a $\rr$-rational real algebraic surface diffeomorphic to $M$.  

By Comessatti classification, if $M$ is orientable its genus must be less than 2. 
It is known that the sphere and the torus each admits a unique rational model modulo biregular isomorphism. Thanks to the latter theorem, there is also only one rational model for the Klein bottle. Hence, the next natural question is:  "is there a genus $h$ for which the nonorientable surface of Euler characteristic $1-h$ admits several rational models?".

\bigskip
One of the main tools used in the proof of Theorems \ref{theo:mainf} and \ref{theo:ruledrat} is a new characterization
of a classical invariant of real algebraic varieties used in the approximation problem in case the target space is the usual sphere.

Given a compact nonsingular real algebraic variety $V$, consider a
smooth projective variety $X$ over $\rr$, such that $V$
and $\xr$ are isomorphic as real algebraic
varieties. We denote by $\halg^2(\xc,\z)$ the subgroup of $H^2(\xc,\z)$ that consists of
the cohomology classes that are Poincar\'e dual to the homology classes in
$H_{2n-2}(\xc,\z)$ represented by divisors in $\xcomp$. We set
 
$$
\hcalg^2(\xr,\z)=i^*(\halg^2(\xc,\z))\;, 
$$
 
\noindent where $i:\xr\hookrightarrow \xc$ is the inclusion map.
We will denote by
$\Gamma (X)$ the quotient group $H^2(\xr,\z)/\hcalg^2(\xr,\z)$. It
is easy to check that for a given nonsingular real algebraic variety $V$, the group $\hcalg^2(\xr,\z)$ does not depend
on the associated variety $X$. We can identify $V$ and $\xr$ and set
$$
\hcalg^2(V,\z)=\hcalg^2(\xr,\z)
$$

We will use the following notation:
$$
\Gamma(V)=\Gamma (X)\;.
$$

There is a close connection between the subgroup
$\hcalg^2(V,\z)$ and the topological closure of the
space $\mathcal {R}(V,S^2)$ in $\mathcal {C}^\infty(V,S^2)$. More
precisely, the following result is well-known, cf. \cite[Chapter
13]{BCR98} and \cite{BBK}.
 
\begin{prop}Let $V$ be a compact nonsingular real algebraic variety.
A given ${\mathcal C}^\infty$ map  $f\colon V\rightarrow S^2$ can be approximated by regular
maps in the ${\mathcal C}^\infty$ topology, if and only if
$f^*(\kappa) \in \hcalg^2(V,\z)$. Here $\kappa$ is a fixed
generator of the group $H^2( S^2,\z)$.
 \end{prop}

In terms of the quotient  $\Gamma$, here are the already known cases:

\begin{theo}[\cite{BK87a,BK87b}] Let $V$ be a $\rr$-rational
real algebraic surface. Then 
$$ 
\Gamma (V)=\begin{cases} \z
&\textrm{ if } V \textrm{ is diffeomorphic to } S^1\times S^1\\
  0&\textrm{ in all other cases.}
\end{cases}
$$
\end{theo}

On the other hand, a smooth projective surface $X$
is a {\em Del Pezzo surface} iff $X$ is
irreducible and the anticanonical divisor $-K_{\xsig}$ is ample.
The degree of $X$ is $\deg(X)=K_{X}^2$.
For Del Pezzo surfaces,  it is known that $1\leq \deg(X)\leq 9$. A real Del Pezzo surface $X$ is $\cc$-rational and is not $\rr$-rational when $\xr$ is not connected.
 
\begin{theo}[\cite{JM03}]\label{theo:dp2}
Let $V$ be a real algebraic surface biregularly isomorphic to the real part of a maximal real Del Pezzo surface of degree~2, then 
$$
\Gamma(V)=\zz\;.
$$
\end{theo}

\begin{theo}[\cite{JM03}]\label{theo:gene}
Let $V=\xr$ be a $\cc$-rational
real algebraic surface. Then 
$$ 
\Gamma (V)=\begin{cases} \z
&\textrm{ if } V\textrm{ is diffeomorphic to } S^1\times S^1\\ 
\zz
&\textrm{ if }
 V \textrm{ is as in
  Theorem~\ref{theo:dp2}} \\ 0&\textrm{ in all other cases.}
\end{cases}
$$
\end{theo}

\medskip
\paragraph{\textbf{Convention.}}
A real algebraic variety is smooth projective and geometrically
irreducible, unless otherwise stated.

\section{Algebraic morphisms to the standard sphere}\label{sec:h2calg}

Let again $i:\xr\hookrightarrow \xc$ be the canonical
injection of the set of real points into the set of complex points of a real algebraic surface.
Consider the induced restriction morphism
$$
i^*\colon H^2(\xc,\z) \to H^2(\xr,\z)\;.
$$
We will use the notation $X^{or}$ for the disjoint union of the orientable connected components of the real part and $X^{nor}$ for the nonorientable part.
The morphism $i^*$ has the natural splitting $H^2(\xc,\z) \to H^2(X^{or},\z)\oplus H^2(X^{nor},\z)$.
Since $X^{nor}$ is nonorientable of dimension 2, $H^2(X^{nor},\z)$ is canonically isomorphic to the group $H^2(X^{nor},\zz)$ by reduction modulo $2$. To see this, apply twice the universal-coefficient theorem \cite[5.5.10]{Spanier}.

We will identify the group $H^2(\xr,\z)$ with the direct sum $H^2(X^{or},\z)\oplus H^2(X^{nor},\zz)$ and still use the notation $i^*$ for the composed morphism 
$$
H^2(\xc,\z) \to H^2(X^{or},\z)\oplus H^2(X^{nor},\zz)\;.
$$
 
The manifolds $\xc$ and $\xr$ are compact and $\xc$ is
orientable.
The Gysin
morphism $i_!$ can be defined by the commutative diagram:
 
\begin{equation}
\begin{CD} H^2(\xc,\z) @>{i^*}>>H^2(X^{or},\z)\oplus H^2(X^{nor},\zz)\\
@V{D_\cc}V{\cong}V @V{D^{or}_\rr\oplus D^{nor}_\rr}V{\cong}V
\\ H_2(\xc,\z) @>{i_!}>> H_0(X^{or},\z)\oplus H_0(X^{nor},\zz)\;,
\end{CD}
\label{diag:gysin}
\end{equation}

where the isomorphisms $D_\cc$, $D^{or}_\rr$ and $D^{nor}_\rr$ come from Poincar\'e duality applied to the orientable 4-dimensional manifold $\xc$, the orientable 2-dimensional manifold $X^{or}$ and the nonorientable 2-dimensional manifold $X^{nor}$.

Let $S$ and $M$ be two transverse oriented submanifolds in an oriented manifold $X$.
We attach $+1$ to a point $P\in S\cap M$ if the orientation of the tangent space $T_P X$ coincide with the orientation given by the direct oriented sum $T_PS\oplus T_PM$ and $-1$ otherwise.
With this convention in mind, we obtain a well-defined class $[S\pitchfork M]$ in $H_0(M,\z)$.
Now if $M$ is nonorientable, the class $[S\pitchfork M]$ is well-defined modulo 2 in $H_0(M,\zz)$ \cite{Hirsch}. The following lemma is an exercise in algebraic topology.

\begin{lem}
\label{lem:trans}
Let $S$ be an oriented 2-dimensional closed
submanifold of $\xc$ transverse to $\xr$, denote by $[S]$ its fundamental class in $H_2(\xc,\z)$, then
$$
i_!([S])=[S\pitchfork X^{or}]\oplus [S\pitchfork X^{nor}]\;.
$$
\end{lem}

Let $X$ be real algebraic surface and suppose that $X(\rr)\ne\emptyset$. Denote by $\{M_j\}_{j\in J}$ the set of connected components of $\xr$. The $\z$-module $H^2(\xr,\z)$ splits into a direct sum $\oplus_{j\in J}H^2(M_j,\z)$. 

Let $J'\subset J$ be a subset and $l$ be a class in $H^2(\xc,\z)$, we will call the image of $l$ by the composed map $H^2(\xc,\z)\to H^2(\xr,\z)\to H^2(\oplus_{j\in J'} M_j,\z)$ the restriction of $i^*(l)$ to $\oplus_{j\in J'} M_j$. 
For a connected component $M_j$ of $\xr$, we will say that a generator class $\eta_j$ of $H^2(M_j,\z)$ belongs to $\hcalg^2(\xr,\z)$ iff the class $1\eta_j\oplus \bigoplus_{a\ne j}0\eta_a$ belongs to $\hcalg^2(\xr,\z)$.

\section{Birational equivalence, orientability and regular maps}
\label{sec:birat}

Let $X$ be a smooth projective surface over $\rr$.
We say that a smooth complex curve $E$ of $\xcomp$ is a {\em {$(-1)$-curve}} if it is rational and $E^2=-1$. If a $(-1)$-curve $E$ is defined over $\rr$, there exists a blowdown $\pi: X\rightarrow Y$ over $\rr$ onto a smooth surface such that $E$ contracts to a real point $P=\pi(E)\in \yr$ (in particular, $\yr$ and $\xr$ are not empty and have the same number of connected components). If $E$ is not defined over $\rr$, then $\sigma (E)$ is another $(-1)$-curve.
If the intersection number $E\cdot\sigma (E)=0$, then we can blowdown over $\rr$ the divisor $E+\sigma (E)$.
We say that a smooth projective surface over $\rr$ is relatively {\em $\rr$-minimal} if it contains neither real $(-1)$-curves nor pairs of disjoint complex conjugated $(-1)$-curves.
If $E\cdot\sigma (E)\ne 0$, then we cannot blowdown $E+\sigma (E)$ over $\rr$ and the surface can be $\rr$-minimal but not $\cc$-minimal, see the next section.

Let $L$ be a {\em real} algebraic curve on an algebraic surface $X$. There are two algebraic bundles naturally associated to $L$. Namely the $\cc$-line bundle $\cE=\ox(L)$ over $\xc$ and the $\rr$-line bundle $\cL$ over $\xr$ satisfying the relation
\begin{equation}
\cL\otimes\cc=\cE_{\vert\yr}\;.
\label{eq:restr}
\end{equation}

We will use  the first Chern class $c_1(L)=c_1(\ox(L))$  in $H^2(\xc,\z)$ and the first Stiefel-Whitney class $w_1(L(\rr))=w_1(\cL)$  in $H^1(\xr,\zz)$.
We denote by $\beta \colon H^1(\xr, \zz)\to H^2(\xr,\z)$ the {\em Bockstein} homomorphism
induced in cohomology by the usual exact sequence
 $$ 
0 \rightarrow \z
\labelto{\times 2}{2} \z \rightarrow \zz \rightarrow 0 \;. 
$$

\begin{lem}
\label{lem:curve}
Let $X$ be a real algebraic surface and let $L$ be a real algebraic curve on $X$, then $i^*(c_1(L))=\beta\circ w_1(L(\rr))$  in $\hcalg^2(\xr,\z)$. In particular, the class $i^*(c_1(L))$ is 2-torsion.
\end{lem}

\begin{proof}
From (\ref{eq:restr}), we get
 $i^*(c_1(\cE))=c_1(\cL\otimes\cc)$ by functoriality of Chern classes.
Since $\cL\oplus\cL$ and $\cL\otimes\cc$ are naturally isomorphic  as oriented real vector bundles and $c_1(\cL\otimes\cc)=\beta\circ w_1(\cL)$ \cite[Problem~15.D and Lemma~14.9]{MiSt}, the lemma follows.
\end{proof}

\begin{cor}
\label{cor:real}
For a real curve $L$, the class $i^*( c_1(L))\ne0$  in $\hcalg^2(\xr,\z)$ if and only if there exists a connected component $M$ of $\xr$ such that $\deg_M( w_1(\cL)^2)$ is odd. In particular, if $L\cdot L$ is odd, $i^*(c_1(L))$  is a nontrivial class of order 2 in $\hcalg^2(\xr,\z)$.
\end{cor}

\begin{proof} The image of the Bokstein homomorphism is given by

$$
\beta(w_1(\cL))=w_1(\cL)\cup w_1(\cL)
$$ 
thanks to the Whitney duality theorem. Furthermore, we have $\deg( w_1(\cL)^2)\equiv L\cdot L \mod 2$. Indeed $\deg( w_1(L(\rr))^2) \equiv L(\rr)\cdot L(\rr) \mod 2$ in $H^1(\xr,\zz)$ and $L\cdot L\equiv L(\rr)\cdot L(\rr) \mod 2$ \cite[Chap. III]{Si89}.
\end{proof}

\begin{rem}
This result was known in case $\xr$ is connected \cite[Th.~12.6.13]{BCR98}.
\end{rem}

\begin{prop}
\label{prop:exceptcurve}
Let $X$ be a real algebraic surface containing a  $(-1)$-curve $L$ defined over $\rr$. There is only one connected component $M$ of $\xr$ meeting $L(\rr)$. Furthermore $M$ must be nonorientable and we have 
$$
\eta_M = i^*(c_1(L))\;.
$$

Hence, the generator class $\eta_M$ of $H^2(M,\z)$ is a nontrivial 2-torsion class of $\hcalg^2(\xr,\z)$.
\end{prop}

\begin{proof}
By Corollary~\ref{cor:real}, $i^*(c_1(L))$  is a nontrivial 2-torsion class.
Another consequence of $L\cdot L=-1$ is that $L$ must have a real point. Let $M$ be a connected component of $\xr$ having a nontrivial intersection with $L(\rr)$. The curve $L$ is smooth and rational, hence have a connected real part. Then $L(\rr)\subset M$ and $M$ must be nonorientable and finally $\eta_M = i^*(c_1(L))$.
 \end{proof}

We recovered the well-known fact that $\Gamma$ is not a birational invariant. But we will prove more: $\Gamma$ is not  even an invariant of relative minimal models, see Theorem~\ref{theo:main} and Corollary~\ref{cor:orconic}.

Given a dominant $\rr$-birational morphism $X\longrightarrow Y$ between smooth real algebraic surfaces, we have a natural injection $\Gamma(X)\hookrightarrow \Gamma(Y)$.
In case $X$ admits a unique minimal model, i.e. when $\kod(X)\geq 0$, we can reduce the computation of $\Gamma(X)$ to that of $\Gamma(Y)$.

In case $X$ is $\cc$-ruled, the result depends on the real minimal model and the dominant map.

\begin{prop}
\label{prop:caneven}
Let $X$ be a real algebraic surface and denote by $\kx$ its canonical line bundle. The class $i^*(c_1(\kx))$ is zero in $\hcalg^2(\xr,\z)$ if and only if the Euler characteristic $\chi(M)$ is even for any nonorientable component $M$.
\end{prop}

\begin{proof}
Indeed, if $\xr\ne\emptyset$, $\kx$ is representable by a real divisor and  $i^*(c_1(\kx))$ belongs to $H^2(X^{nor},\z)$ by Lemma~\ref{lem:curve}. We have then
$$
i^*(c_1(\kx))\equiv _{\mathrm {mod}\ 2} w_2(T\xr)
$$
and we conclude by \cite[Cor.~11.12]{MiSt} about Stiefel-Whitney numbers.
\end{proof}

\section{Real conic bundles over curves}
\label{sec:conic}

Let $X$ be a smooth real algebraic surface and $B$ a smooth real algebraic curve.
A connected component $M$ of $\xr$ is said to be a {\em spherical} (resp. {\em torus}, resp. {\em Klein}) component if $M$ is diffeomorphic to the sphere $S^2$ (resp. the torus, resp. the Klein bottle).

\begin{defin}
A morphism  $\rho \colon X\to B$ is  a {\em ruling} iff the generic fiber is isomorphic to $\pp^1$.
The morphism $\rho$ is a
{\em conic bundle} iff every fiber is isomorphic to a plane conic.
\end{defin}

When the map $\rho$ is defined over $\rr$, we will say that a fiber of $\rho$ is {\em real} if it is located over $B(\rr)$ and {\em imaginary} otherwise.

A  ruling is {\em $\cc$-minimal} iff no fiber contains a  $(-1)$-curve. A  real ruling is {\em $\rr$-minimal} iff no real fiber contains a  {\em real} $(-1)$-curve and no imaginary fiber contains a $(-1)$-curve.
A $\cc$-minimal real ruling is clearly $\rr$-minimal but the converse does not hold in general.
A $\cc$-minimal ruling is isomorphic to a locally trivial $\pp^1$-bundle.
A $\rr$-minimal ruling is a real conic bundle. 

Recalling that $\kk$-ruled means $\kk$-birational to a product $\pp^1\times B$, we have:
\begin{prop}
A given $\rr$-minimal, $\cc$-ruled and non $\cc$-rational real algebraic surface  is $\rr$-ruled if and only if it is $\cc$-minimal.
\end{prop}

A surface $X$ endowed with a minimal ruling (over $\cc$ or over $\rr$) in the above sense is not a minimal model in the sense of Mori theory, it is only relatively minimal when $g(B)\ne 0$. 
Indeed, there exist many birationally equivalent minimal rulings. 
A birational equivalence between two minimal ruling is a composition of elementary transformations. Over $\cc$, an elementary transformation centered at a point $p$ is the blow-up centered at $p$ composed by the contraction of the strict transform of the fiber containing $p$.

Over $\rr$, there are two kinds of elementary transformations. We denote by $\elm_p$ the elementary transformation centered at a {\em real} point $p$ that is smooth in $X_{\rho(p)}$ and by $\elm_{p,\sigma(p)}$ the composition of the elementary transformations centered at  $p$ and $\sigma(p)$ provided that $X_{\rho(p)}$ and $X_{\rho(\sigma(p))}$ are distinct conjugated fibres.

We will use the following classification theorem:

\begin{theo}
\label{theo:classif}
Let $X$ be a smooth real surface with a real ruling $\rho \colon X\to B$ over a smooth real curve $B$.
Then $X$ is $\rr$-birational to the smooth $\rr$-minimal projective completion $X^g$ of the real conic bundle defined in some affine open subset of $\bb{A}^2\times B$ by an equation
\begin{equation}
x^2+y^2=g(z)\label{eq:conic}\;,
\end{equation}
where $g$ is a real rational function over $B$ with no pole in $B(\rr)$, and whose all real zeros are simple.
\end{theo}

\begin{proof}
See \cite[V.2]{Si89} and \cite[VI.3]{Si89}.
\end{proof}

\begin{rem}Note the number of real zeros of $g$ that belong to a connected component $B_1$ of $B(\rr)$ is even.
Indeed, the function $g$ changes sign in the neighborhood of a zero in the topological circle $B_1$.
\end{rem}

\begin{prop}
Denote by $n$ the number of connected components of $B(\rr)$ and by $2s$ the number of real zeros of $g$. Then the real part $X^g(\rr)$ of $X^g$ is diffeomorphic to the disjoint union of $t$ tori and $s$ spheres, where $t$ satisfies $t\leq n$.
\end{prop}

\begin{proof}
From Equation~\ref{eq:conic}, the topology of $X^g(\rr)$ is easy to understand. The real
zeros $\{z_l\}_{1\leq l\leq 2s}$ of the function $g$ determine $s$
connected arcs in $B(\rr)$ over which the real fibers
$X^g_z(\rr)$ are not empty. Over each of these arcs, there
is a connected component of $X^g(\rr)$ which is homeomorphic to a
sphere. The torus components are located over components of $B(\rr)$ where $g$ is strictly positive.
\end{proof}
 
Let $M$ be a spherical component of $X^g(\rr)$ and $X^g_{z_l}$ be a real fiber over a zero of $g$ such that $X^g_{z_l}(\rr)\subset M$. The real singular fiber $X^g_{z_l}$ is the union of two complex
conjugated $(-1)$-curves $E$ and $\sigma(E)$ whose intersection point $p$
is the only real point of the fiber.
 
\begin{lem}
\label{lem:nsalg}
Let $\xsig$ be a $\cc$-ruled surface defined over $\rr$, then
$$
\hcalg^2(\xr,\z)=\im i^*\;.
$$
\end{lem}

\begin{proof}
By definition, we have $\hcalg^2(\xr,\z)= i^*(\ns(\xcomp))$.
Moreover, for a complex nonsingular variety $V$, we have the long
exact sequence coming from the exponential exact sequence.
In addition, considering the isomorphism $\pic(V)\cong
H^1(V,\mathcal{O}^*)$, we obtain the following exact sequence
 
$$\cdots\to H^1(V,{\mathcal O})\to
\pic(V)\labelto{c_1}{1}H^2(V,\mathbb Z)\to H^2(V,{\mathcal
O})\to\cdots$$
 
\noindent The Lemma is now clear since $\ns(V)=c_1(\pic(V))$ and, for a $\cc$-ruled surface $V$, $\dim H^2(V,\mathcal{O})=0$.
\end{proof}

Let $X$ be a $\rr$-minimal conic bundle, by Theorem~\ref{theo:classif}, there exist a surface $X^g$ and a finite sequence $T$ of real elementary transformations such that $T(X^g)=X$. 

The  connected components $\{M_j\}_{j\in J}$ of $\xr$ are then spherical, toral or of Klein type.
 We call respectively $S\subset J,T\subset J,K\subset J$ the subsets of indexes corresponding to the spherical, torus and Klein components respectively. In particular, we have $X^{nor}=\oplus_{j\in K}M_j$. 

Now, for each spherical component $M_j$ of $\xr$,  there exist two singular fibers $E_j^c+\sigma (E_j^c)$, $c=1,2$, such that $E_j^c\cap \sigma (E_j^c)\in M_j$ is reduced to a single real point. Let us denote by $\cE_j^c$ the associated algebraic $\cc$-line bundles and by $N$ the submodule of $H^2(\xc,\z)$ generated by the $2s$ classes $\{ c_1(\cE_j^c)\}_{c\in\{1,2\},j\in S}$.

\begin{lem}
\label{lem:ns}
Let $X$ be a $\rr$-minimal conic bundle. 
\begin{enumerate}
\item The N\'eron-Severi group $\ns(\xcomp)$ is generated by $N$, the class $f$ of a fiber and the class $h$ of a section. 
\item The canonical class is given by
$$
c_1(\mathcal{K}_{X}) =  rf-2h+\sum_{c\in\{1,2\},j\in S} c_1(\cE_j^c)
$$
for some integer $r\in \z$.

\end{enumerate}
\end{lem}

\begin{proof}
Indeed, over $\cc$, we can blow down the curves $E_j^c$ to obtain a $\cc$-minimal complex ruled surface $\rho' \colon Y\longrightarrow B$ whose N\'eron-Severi group is generated by the class $f'$ of a fiber  and the class $h'$ of a section. Moreover, the canonical class $c_1(\mathcal{K}_{Y})$ is a linear combination 
\begin{equation}\label{eq:can}
rf'-2h'
\end{equation}
for some integer $r\in \z$ \cite[III.18]{Be}.

The  strict transform of a generic fiber of $Y$ is a fiber of $X$ and the strict transform of a section of $Y$ is a section of $X$.
Furthermore, for a blow-up centered at $p\in Y$, the total transform of the fiber $Y_{\rho'(p)}$  is a singular fiber for $X\to B$ of the form $E+\sigma (E)$, where $E$ is a $(-1)$-curve. 

Hence, the group $\ns(\xcomp)$ is generated by the classes $c_1(E_j^c)$, the class $f$ of a fiber and the class $h$ of a section.
Furthermore, we deduce from (\ref{eq:can}) that $c_1(\mathcal{K}_{X}) =  rf-2h+\sum c_1(\cE_j^c)$.
\end{proof}

\begin{lem}
\label{lem:fibre}
Given any $\rr$-minimal conic bundle $X\to B$, the canonical class and the class of a fiber satisfy $i^*(\kx)=0$ and  $i^*(f)=0$ in $H^2(\xr,\z)$.
\end{lem}

\begin{proof}
We may assume that $\xr\ne\emptyset$ hence there exist a fiber $F$ of $\rho$ and a canonical divisor which are real. Then $F.F=0$ as a fiber and the conclusion about $f$ follows from Corollary~\ref{cor:real}. 
By Lemma~\ref{lem:curve}, $i^*(\kx)$ is a 2-torsion  class, hence trivial when restricted to an orientable component. Furthermore, the restriction of $i^*(\kx)$ to a Klein component $M_j, j\in K$, is trivial as $w_2(M_j)=0$.
 \end{proof}

\begin{lem}
\label{lem:n}
For each spherical component $M_j\subset \xr$, any generator class $\eta_j$ of $H^2(M_j,\z)$ belongs to $\hcalg^2(\xr,\z)$. More precisely, we have 
$$
i^*(N) = \oplus_{j\in S}H^2(M_j,\z)\;.
$$
\end{lem}

\begin{proof}
By Theorem~\ref{theo:classif}, there exist a real ruling $\rho' \colon X^g\longrightarrow B$ and a finite sequence of real elementary transformations $T(X^g)=X$, where $X^g$ is the $\rr$-minimal projective completion of the conic bundle defined by
$$
\{(x,y,z)\in \bb{A}^2\times B \,\vert\, x^2+y^2=g(z)\}\quad\textrm{ and }\quad \rho' (x,y,z)=z\;.
$$

As $\rho$ and $\rho'$ are $\rr$-minimal, we may assume that there is no center of elementary transformation of $T$ that belongs to a reducible fiber. In particular, restricted to a neighborhood of a spherical component $M$ of $X^g(\rr)$ in $X^g(\cc)$, $T$ is a real isomorphism.

Let $X^g_z$ be a real fiber over a zero of $g$ such that $X^g_z(\rr)\subset M$. The real singular fiber $X^g_z$ is the union of two complex
conjugated lines $E$ and $\sigma(E)$ whose intersection point $p$
is the only real point of the fiber. The tangent plane to $X^g(\rr)$ at
$p$ is generated by $\frac{\partial}{\partial x_1}$ and
$\frac{\partial}{\partial y_1}$, where $x=x_1+ix_2$ and
$y=y_1+iy_2$. It is easy to check that the tangent plane to $E$ at
$p$ is generated by $i\frac{\partial}{\partial
x_1}-\frac{\partial}{\partial x_2}$ and $\frac{\partial}{\partial
y_1}+i\frac{\partial}{\partial y_2}$.
Then $E$ is transverse to $X^g(\rr)$ at $p$ in $X^g(\cc)$. 

Hence $T(E)$ is transverse to $X(\rr)$ at $T(p)$ in $X(\cc)$.
By
Lemma~\ref{lem:trans}, the image of $[T(E)]$ by the Gysin morphism $i_!$ is a generator of 
$$
H_0(T(M),\z)\hookrightarrow H_0(X^{or},\z)\oplus H_0(X^{nor},\zz)\;.
$$

Therefore, we conclude by using the commutative diagram \ref{diag:gysin}.
\end{proof}

\begin{lem}
\label{lem:torus}
Let $M_j$ be a torus component. Then for any generator class $\eta_j$ of $H^2(M_j,\z)$ we have $\eta_j \not\in\hcalg^2(\xr,\z)$.
\end{lem}

\begin{proof}
By Lemmas~\ref{lem:ns}, \ref{lem:fibre} and \ref{lem:n}, it suffices to prove that the restriction of $i^*(h)$ to a torus component is trivial.

The canonical class $c_1(\kx)$ is a linear combination $rf-2h+\sum  c_1(\cE_j^c)$ for some integer $r\in \z$. From Lemma~\ref{lem:n} and Lemma~\ref{lem:fibre}, the restrictions of $i^*(2h)$ and $i^*(\kx)$ to a torus component are equal.
Moreover, the class $i^*(\kx)$ is trivial when restricted to an orientable component. The restriction of $i^*(h)$ to a torus component is then a 2-torsion class, hence trivial.
\end{proof}

\begin{lem}
\label{lem:section}
Let $h$ be the class of a section on a $\rr$-minimal conic bundle $X$, then the restriction of $i^*(h)$ to the nonorientable part is the class $\oplus_{j\in K} \eta_j\in H^2(\xr,\z)$.
\end{lem}

\begin{proof}
The restriction of $i^*(h)$ to the nonorientable part is 2-torsion. We will prove that we can choose a section $H$ transverse to $\xr$ such that $\#(H\cap M)$ is odd for any nonorientable component. The conclusion will then follow from Lemma~\ref{lem:trans}.

As in the proof of lemma~\ref{lem:n}, we will use the surface $X^g$ and the transform $T$. Let $\Sigma$ be the finite set of real centers of elementary transformations of $T$. If $M$ is spherical, $T(M)$ is also spherical. If $M$ is a torus component, then $T(M)$ is a torus component when $\#(\Sigma \cap M)$ is even and a Klein component if $\#(\Sigma \cap M)$ is odd.

Let $H'$ be a section of $X^g$. The curve $H=T(H')$ is then a section of $X$. Since $\Sigma$ is finite, we can move $H'$ to ensure that for all $z\in\rho'(\Sigma)$, $X^g_z(\rr)\cap H'(\cc)=\emptyset$. Hence the point $p=T(X^g_z)$ is real and belongs to the intersection $H\cap \sigma(H)$.

For each $z\in\rho'(\Sigma)$, the intersection $H\cap \sigma(H)$ is transverse at $p$ and real. Hence $H$ is transverse to $\xr$ at $p$. If necessary, we can perturb $H$ to obtain transversality to $\xr$ at each point. 

Now for a non spherical component $T(M)$ of $X$, the degree of the restriction of $i_!([H])$ to  $T(M)$ is equal to the sum of the degree of $i_!([H'])_{\vert M}$ and $\#(\Sigma \cap M)$. 

By Lemma~\ref{lem:torus}, $i_!([H'])_{\vert M}=0$ for a torus component $M$ of $X^g(\rr)$ and $\#(\Sigma \cap M)$ is odd when $T(M)$ is a Klein component of $\xr$. The conclusion follows.
\end{proof}

\medskip
Given a $\rr$-minimal conic bundle $X$, from Lemmas~\ref{lem:nsalg} to \ref{lem:section} we get
$$
\hcalg^2(\xr,\z)=\left<\bigoplus_{j\in K} \eta_j, \{\eta_j;j\in S\}\right>\;.
$$

In other words, the group $\hcalg^2(\xr,\z)$ is generated by the spherical classes and the sum of all the Klein classes. We deduce the theorem:

\begin{theo}
\label{theo:main}
Let $X\to B$ be a $\rr$-minimal conic bundle.
Denote by $t$ the number of torus components of $\xr$ and by $k$ the number of Klein components.  Then
$$
\Gamma (X)=\z^t\oplus(\zz)^{k-1}\;.
$$
\end{theo}

\begin{cor}
\label{cor:orconic}
Given a real ruling $X\to B$ on a surface with orientable real part $\xr$, we have $\Gamma(X)=\z^t$, where $t$ is the number of torus components.
\end{cor}

\begin{proof}
The orientability of the real part implies that the  real ruling gives rise to a $\rr$-minimal conic bundle by making contractions in imaginary fibers only.
\end{proof}

\section{Surfaces of negative Kodaira dimension}\label{sec:negkod}

To prove Theorems \ref{theo:mainf} and \ref{theo:ruledrat}, we will use the following (recall that the subgroup $\halg^1(V,\zz)\subset H^1(V,\zz)$ is generated by the cohomology classes Poincar\'e dual to the homology classes represented by Zariski closed algebraic hypersurfaces of $V$):

\begin{theo}[\cite{Ku}]
\label{theo:ku}
Let $V$ be a compact nonsingular real algebraic variety and $W$ be a compact connected nonsingular rational real algebraic surface.
Given a ${\mathcal C}^\infty$ map  $f\colon V\rightarrow W$, the following conditions are equivalent:

\begin{enumerate}
\item $f$ can be approximated by regular maps;
\item $f$ is homotopic to a regular map;
\item either $W$ is diffeomorphic to a sphere and 
$$
f^*(H^2(W,\z))\subset \hcalg^2(V,\zz)\;,
$$
or $W$ is not diffeomorphic to a sphere and 
$$
f^*(H^1(W,\z))\subset \halg^1(V,\zz)\;.
$$
\end{enumerate}
\end{theo}
 
Given a real algebraic surface $X$, we denote by $t$ the number of torus components, and by $k$ the number of Klein components. In case $X$ admits a real ruling $\rho  \colon X \to B$, we denote by $k'$ the number of Klein components of $\xr$ whose image by $\rho$ is a connected component of $B(\rr)$.

 \begin{theo}
\label{theo:maingamma}
If $X$ is a $\cc$-ruled non $\cc$-rational real algebraic surface, then
$$
\Gamma (X)=\z^t\oplus(\zz)^{k'-1}\;.
$$ 
 \end{theo}
 
\begin{proof}
Let $Y$ be a $\rr$-minimal model of $X$, thus the number of Klein components of $\yr$ is exactly $k'$ and by Theorem~\ref{theo:main}, $\Gamma(Y)=\z^t\oplus(\zz)^{k'-1}$. The conclusion follows from Proposition~\ref{prop:exceptcurve}.
\end{proof}

\begin{theo}\label{theo:uniruledor}
Let $V=\xr$ be an orientable  real algebraic surface with $\kod(X)=-\infty$. Denoting by $t$ the number of components diffeomorphic to a torus,  we have 
$$
\Gamma(X)=\z^t\;,
$$ 
except in case $X$  is the maximal real Del Pezzo surface of degree 2 for which $\xr$ is the disjoint union of 4 spheres and $\Gamma(X)=\zz$.
\end{theo}

\begin{rem}
It is an amazing fact that the torus components measure the obstruction to approximate differentiable maps. Indeed, the only case that is known so far is  the rational torus $S^1\times S^1$ realized as the real part of the quadric surface $\pp^1\times\pp^1$ endowed with the usual real structure. The proof of  $\Gamma(X)=\z$ uses the torus decomposition as a product of real algebraic curves.  
\end{rem}

\begin{proof}
Since $V$ is orientable, we may assume that $X$ is $\rr$-minimal. Thus, we get the conclusion from Theorem~\ref{cor:orconic} in case $X$ admits a real ruling and from theorems~\ref{theo:dp2} and \ref{theo:gene} when $X$ is $\cc$-rational.
\end{proof}

\begin{cor}
Let $V=\xr$ be an orientable surface  of negative Kodaira dimension which is not biregular to a maximal real Del Pezzo surface of degree 2. 

The space of regular maps ${\mathcal R}(V,{\mathbb S}^2)$ is
dense in the space of ${\mathcal C}^\infty$maps ${\mathcal
C}^\infty(V,{\mathbb S}^2)$ if and only if all the connected components of $V$ are spherical.
\end{cor} 

\begin{proof}[Proof of Theorem~\ref{theo:mainf}]
Theorem~\ref{theo:mainf} follows from Theorem~\ref{theo:ku} and Theorem~\ref{theo:maingamma}
\end{proof}

\begin{proof}[Proof of Theorem~\ref{theo:ruledrat}]
We got the conclusion in case $W$ is diffeomorphic to a sphere from Theorem~\ref{theo:mainf} and when $W$ is not diffeomorphic to a sphere by Theorem~\ref{theo:ku}, and the fact that a connected $\cc$-ruled real algebraic surface $V$ satisfies  $\halg^1(V,\zz)=H^1(V,\z)$  \cite{Ab00,M03}.
\end{proof}

\section{Rational Klein bottles}\label{sec:klein}

This short section is devoted to the proof of Theorem~\ref{theo:klein}.

\begin{theo}
The Klein bottle admits a unique rational model. Namely, the real part of the real Hirzebruch surface $F(1)$.
\end{theo}

Here we can use indifferently  the words $\cc$-rational or rational because connected $\cc$-rational surfaces are $\rr$-rational \cite{Si89}.

\begin{proof} 
We want to prove that $M$ is biregularly isomorphic to the real part of the Hirzebruch surface $F(1)$.
Let $M$ be a $\cc$-rational real algebraic surface diffeomorphic to the Klein bottle. Let $X$ be a $\rr$-minimal smooth projective complexification of $M$. As $M\cong X(\rr)$ is connected,  $X$ is $\cc$-minimal and it is a Hirzebruch surface $F(n)$. Furthermore,  $n$ is odd and $n>1$. Indeed, the only $\cc$-minimal $\rr$-rational surfaces are the real Hirzebruch surfaces $F(n)$ with $n\ne 1$ and $F(n)(\rr)$ is nonorientable if and only if $n\equiv 1 \mod 2$ \cite{Si89}.

Let us denote by $H$ the unique section of the natural real ruling $\rho \colon X\to \pp^1$  such that $H^2=-n$. Choose $\frac {n-1}2$ points $p_1,\dots,p_{\frac {n-1}2}$ of $H$ that belong to imaginary fibres of $\rho$ and let $X'= \elm_{p_1,\sigma(p_1)}\circ \cdots \circ \elm_{p_{\frac {n-1}2},\sigma(p_{\frac {n-1}2})}(X)$.

Then $X'(\rr)$ is biregularly isomorphic to $X(\rr)$ and $n'=n-2(\frac {n-1}2)$.
Furthermore, the transformed surface $X'$ of $X$ is $\cc$-isomorphic to $F(1)$.
\end{proof}


\end{document}